\documentclass[11pt,english]{article}

\usepackage{graphicx,psfrag,epsfig,amsfonts,amssymb,amsmath,babel}

 \newtheorem{theorem}{Theorem}[section]
\newtheorem{corollary}[theorem]{Corollary}
\newtheorem{lemma}[theorem]{Lemma}
\newtheorem{proposition}[theorem]{Proposition}

\newtheorem{definition}[theorem]{Definition}

\newtheorem{prueba}[theorem]{}

\DeclareMathOperator{\dist}{dist}

\begin{document}

\title{Equilibrium States and SRB-like measures of $C^1$ Expanding Maps of the Circle.}

\author{Eleonora Catsigeras\thanks{Instituto de Matem\'{a}tica y Estad\'{\i}stica Rafael Laguardia (IMERL),
 Fac. Ingenier\'{\i}a,  Universidad de la Rep\'{u}blica,  Uruguay.
 E-mail: eleonora@fing.edu.uy
 Address: Herrera
  y Reissig 565. Montevideo. Uruguay. EC was partially supported by CSIC of Universidad de la Rep\'{u}blica and ANII of Uruguay.} \  and Heber Enrich}

\date{}
\maketitle

\begin{abstract}
For any $C^1$ expanding map $f$ of the circle we study the  equilibrium states for the potential $\psi=-\log |f'|$. We formulate  a   $C^1$ generalization of Pesin's Entropy Formula  that holds for all the SRB measures if they exist, and for all the (necessarily existing) SRB-like measures. In the $C^1$-generic case  Pesin's Entropy Formula holds for a unique  SRB measure which is not absolutely continuous with respect to Lebesgue. The result also stands in  the non generic case for which no SRB measure exists.
\end{abstract}

\small

 \noindent {\em Math. Subj. Class. \em (2010): }
Primary 37A05; Secondary 28D05.


\normalsize

\section{Introduction}

For any map $f$ on a compact manifold, if no invariant measure is equivalent to Lebesgue, or if   $f$ is non ergodic with respect to    the invariant  measures that are equivalent to Lebesgue, many substitutive concepts of natural invariant measures have been defined. They translate the statistical asymptotic behavior of Lebesgue-positive sets of orbits, into spatial probabilities. Nevertheless, except under specific conditions in the $C^{1 + \alpha}$ scenario, those statistically good measures do not necessarily coincide, and moreover, they do not necessarily exist. For instance in \cite{blankbunimovich} and \cite{misiurewicz}  the \em natural measures \em are defined as the weak$^*$ limit (if it exists) of the averages $(1/n) \sum_{j= 0}^{n-1} (f^*)^j \nu$ for any probability $\nu \ll m$, where $m$ is the Lebesgue measure and $f^*$ denotes the pull back operator in the space of Borel-probabilities. Similarly, \em SRB measures \em  are defined as the weak$^*$ limit (if it exists) of the averages $\sigma_n(x):= (1/n) \sum_{j= 0}^{n-1} (f^*)^j \delta_x$ for a Lebesgue-positive set of initial states $x$. In \cite{misiurewicz}   a method is exhibited to construct $C^0$ non singular expanding maps of the circle $S^1$, for which there exists a unique natural measure with respect to $m$, and nevertheless, for $m$-almost all the points $x \in S^1$ the averages $\sigma_n(x)$ are non convergent. Thus, even in the case of topological expansion, the notions of SRB and natural measures are different. In \cite{keller2}, \cite{zweimuller}, \cite{jarvenpaatolonen} diverse maps are constructed without any natural limit measure, but with other good ergodic properties (for instance, the existing of a mixing probability). In a general context, neither the existence of natural measures nor of SRB measures   is required for a map $f$    exhibit  statistically good properties with respect to the Lebesgue measure (\cite{CE}).

A third notion of good measure  from the statistical viewpoint  raises from the thermodynamic formalism when considering, if it exists, a probability $\mu$ for which Pesin's Entropy Formula holds \cite{pesin}.
If the hypothesis of   $C^2$   (or   $C^{1 + \alpha}$) regularity is added,   plenty of tight relations were proved among the SRB measures, the absolute continuity with respect to Lebesgue (of the conditional measures along the unstable manifolds), and Pesin's Entropy Formula. See for instance \cite{ledrappieryoung}, \cite{qianzhu}, \cite{suntian}, \cite{araujo}, \cite{araujotahzibi}. But to prove those results, the $C^1$-plus H\"{o}lder regularity is essential.  In the $C^1$ scenario,  generic volume preserving diffeomorphisms still have an invariant measure satisfying Pesin's Entropy Formula \cite{tahzibi}, \cite{bochiviana}. But contrarily to the situation of the $C^{1 + \alpha}$  maps, the $C^1$-generic dynamical systems, under  some hyperbolic-like assumptions, have no invariant measure $\mu$ being (either $\mu$  or its conditional unstable measures) absolutely continuous with respect to Lebesgue \cite{avilabochi}, \cite{bruinhawkins}. Nevertheless, those $C^1$-generic systems still have a unique SRB measure \cite{cq},  \cite{qiuhao}.

Along this paper we consider the family ${\mathcal E}^1$ of all the $C^1$ expanding maps of the circle $S^1$. We recall that a $C^1$ map $f\colon  S^1 \mapsto S^1$ is expanding if $|f'(x)| >1$ for all $x \in S^1$. We denote $ {\mathcal E}^{1 + \alpha} = {\mathcal E}^1 \cap {\mathcal C}^{1 + \alpha}$. Namely $f \in {\mathcal E}^{1 + \alpha}$ if and only if $f \in {\mathcal E}^1$  and besides $f'$ is $\alpha$- H\"{o}lder continuous. We will focus     on the  systems in ${\mathcal E}^1 \setminus {\mathcal E}^{1 + \alpha}$. Our purpose is to state and prove a reformulation of Pesin's Entropy Formula  for these systems, including the non generic ones for which no SRB exists.

For $C^1$ systems,  Ruelle's Inequality \cite{ruelleIneq}  states that for any $f$-invariant   probability measure $\mu$ on the Borel $\sigma$-algebra of $S^1$, the corresponding measure theoretic entropy $h_{\mu}(f)$  satisfies:


\begin{equation}\label{equationRuelleInequality} h_{\mu}(f) \leq \int \log |f'|\, d \mu.\end{equation}


 \noindent Therefore,      $h_{\mu}(f) - \int \log|f'|\, d \mu \leq 0$. By definition, Pesin's Entropy Formula holds   if  the latter difference is equal to zero:


  \begin{equation}
 \label{equationPesinsFormula}
 h_{\mu}(f) = \int \log|f'|\, d \mu \end{equation}


   For any map $f \in {\mathcal E}^{1 + \alpha}$,   \cite{pesin} and   \cite{ledrappieryoung}  prove that  Formula (\ref{equationPesinsFormula})    holds  if and only if   $\mu \ll m$, where $m$ is the Lebesgue measure.   On the contrary, as said above, if    $f$ is only $C^{1}$  then  ${\mathcal E}^1$-generically $f$ has no invariant measure  $\mu$ such that $\mu \ll m$ \cite{cq}. Two questions arise: First, do there exist, for any $f \in {\mathcal E}^1$,   invariant probability measures satisfying   Formula (\ref{equationPesinsFormula})? From the thermodynamic formalism, the answer to this question is known to be affirmative, since $f$ is topologically expansive.  Second, what statistical properties do those probabilities exhibit with respect to the (non invariant) Lebesgue measure? In Theorem \ref{maintheorem} of this paper    we give a   statistical simple description of a nonempty subset of invariant measures that satisfy Formula (\ref{equationPesinsFormula}). We call that description the SRB-like property  \cite{CE}. As a Corollary, if the measure that satisfies Formula (\ref{equationPesinsFormula}) is unique, then it is SRB, also when ${\mathcal E}^1$-generically it is mutually singular with respect to the Lebesgue measure. Besides, if there are physical measures, all of them satisfy Formula (\ref{equationPesinsFormula}). Finally, if SRB measures do not exist, there still exist (uncountably many) probability measures that are distinguished from the general invariant measures by  a weak physical condition, which is similar to the statistical property of SRB measures, and that besides  satisfy Pesin's Entropy Formula.

Even if we conjecture that the results are also true for $C^1$ expanding maps in any dimension, the proofs along this paper work only on one-dimensional compact manifolds. In fact, in Lemma \ref{lemmaEntropy}, we use that there exists a partition of the ambient manifold whose pieces have arbitrarily small diameters, and such that the measure of the union of the boundaries of its pieces is zero  for all the invariant probability measures. This property is trivially satisfied by any   one-dimensional map  whose set of periodic orbits is, at most, countable.

\section{Definitions and Statement  of the Result.} \label{sectionStatements}

    The classic thermodynamic formalism defines the pressure $P_f$ with respect to the potential $$ \psi:= -\log |f'|$$    by


$$P_f = \sup_{\mu \in {\mathcal M}_f} \{ h_{\mu}(f) - \int \log|f'|\, d \mu\},  $$


\noindent where ${\mathcal M}_f$ is the set of all the $f$-invariant Borel probabilities in $S^1$. For  any $f \in {\mathcal E}^1 $ the pressure $P_f$ is   equal to zero (see  \cite{Q3}).
Let us denote with $ES_f$ the (a priori maybe empty) set of all the $f$-invariant probability measures $\mu $ that realize the pressure  $P_f$ as a maximum equal to zero. Precisely:


\begin{equation}\label{equationRLYequality}\mu \in ES_f \ \ \ \mbox{ if and only if } \ \ \ h_{\mu}(f) = \int \log |f'| \, d \mu .\end{equation}


\noindent Namely, the set $ES_f$ is the set of invariant measures that satisfy Pesin's Entropy Formula (\ref{equationPesinsFormula}) of the entropy. The thermodynamic formalism (see for instance   \cite{ke}) for expansive maps states  that $ES_f$ is weak$^*$ compact and convex in the space ${\mathcal M}$ of all the Borel probabilities in $S^1$,  and, if nonempty,  its extremal points are ergodic measures. The  measures in $ES_f$ are called \em equilibrium states \em of $f$ for the $C^0$ potential $\psi= -\log|f'|$.

 Let us recall  some definitions from the statistical  viewpoint.  Consider  for each initial point $x \in S^1$, the following sequence of measures   $\{\sigma_n(x)\}_{n \geq 1}$, that are called  \em empirical probabilities. \em In general they are  non $f$-invariant:


\begin{equation}
\label{eq1}
\sigma_n(x) := \frac{1}{n} \sum_{j= 0}^{n-1} \delta_{f^j(x)}.\end{equation}


\noindent In the   above definition $\delta_y$ denotes the Dirac-delta measure supported on   $y $.
\begin{definition}
 \label{definitionPhysical} \em
 We call a Borel probability measure $\mu$ \em physical or SRB \em if

 \begin{equation}
 \label{equationBasinOfAttraction}B(\mu):= \{ x \in S^1: \ \lim_{n \rightarrow +  \infty} \sigma_n(x) = \mu\} \end{equation}


\noindent has positive Lebesgue measure.  (In the definition   of the set $B(\mu)$  the limit of the measures is taken in the space ${\mathcal M}$ of all the probability measures, endowed with the weak$^*$ topology.)

We call  $B(\mu)$ the \em basin of attraction \em of the physical measure  $\mu$.

\end{definition}

It is standard to check that any physical measure is   $f$-invariant. After the definition above, if there exist physical measures, then they describe  the spatial probabilistical distribution in $S^1$ of the asymptotic behavior of the empirical distributions in Equality (\ref{eq1}),  for a Lebesgue-positive set $B(\mu) \subset S^1$ of initial states. This is the  \em physical \em role of the SRB measures from the statistical viewpoint.
As said in the introduction, the existence and uniqueness of an SRB measure  $\mu$ are  generic properties for   $f \in {\mathcal E}^1$, but $\mu$ is mutually singular with respect to the Lebesgue measure \cite{cq}.  On the other hand, for any  $f \in {\mathcal E}  ^{1 + \alpha}$, the existence and uniqueness of the SRB measure $\mu$ is a well established fact (Ruelle's Theorem). Besides, in this case $\mu$ is equivalent  to Lebesgue and it is the unique equilibrium state for the potential $\psi= - \log |f'|$. Namely, it is the unique probability that satisfies Pesin's Entropy Formula (\ref{equationPesinsFormula}).
In Theorem \ref{maintheorem}  we prove a   generalization of Ruelle's Theorem and of Pesin's Entropy Formula (\ref{equationPesinsFormula})  for any $C^1$ expanding map of the circle. We apply the definition of SRB-like measure, instead of considering only  SRB measures. The gain in this   generalization is that the SRB-like measures always exist. Besides, they still preserve a physical-like meaning (see Proposition \ref{propositionPhysicalLike}) as SRB measures do,  and also, they are   equilibrium states for $-\log|f'|$,  regardless whether SRB measures exist  and whether   such an equilibrium state  is unique.

Before stating   the  precise result,  we need to   revisit the definition of SRB-like  measure. In brief, the nonempty set ${\mathcal O}_f$ of the \em SRB-like probability measures  \em (defined for any continuous map acting on a compact manifold) is the minimal weak$^*$-compact nonempty set of ${\mathcal M}$ that contains all the limits of the convergent subsequences of (\ref{eq1})  for Lebesgue-almost all the initial states $x \in S^1$ (see Definition \ref{definitionobservable}).  Immediately, if there exist SRB measures,   they are SRB-like; and there exists a unique SRB-like measure if and only if there exists a unique SRB probability $\mu$ and   its basin $B(\mu)$ has full-Lebesgue measure. But besides, in the cases that no SRB measure exists, the SRB-like measures still exist and preserve  the  statistical role that the nonexisting SRB measures would exhibit.  Although the construction above is global, each SRB-like measure $\mu \in {\mathcal O}_f$   preserves an individual  weakly physical  meaning, independently of the other measures in the set ${\mathcal O}_f$. This is stated in the following Proposition \ref{propositionPhysicalLike}.
It gives a  characterization of   the SRB-like measures.
To state Proposition \ref{proofPropositionPhysicalLike}, and to argue along the paper,  the space    ${\mathcal M}$ of all the Borel probabilities on $S^1$ is endowed with the weak$^*$ topology.
For each point $x \in S^1$ we denote:


 \begin{equation}\label{equationpomegax}p\omega(x) = \{\mu \in {\mathcal M}: \ \exists \ n_i \rightarrow + \infty \mbox{ such that } \lim _{i \rightarrow + \infty} \sigma_{n_i}(x) = \mu \} \end{equation}


 \noindent where $\sigma_{n}(x)$ is the empirical probability defined in Equality (\ref{eq1}). The set $p\omega(x)$ is the limit set in   ${\mathcal M}$   of the empirical sequence with initial state $x$. We call $p \omega (x)$ \em    the p-limit set of $x$.\em
We fix any weak$^*$-metric in   ${\mathcal M}$. We  denote    by $\mbox{dist}$ this metric.


  \begin{proposition} \label{propositionPhysicalLike}
  A probability measure $\mu$ is SRB-like if and only if for all $\epsilon >0$ the following set $A_{\epsilon}(\mu) \subset S^1$(called basin of $\epsilon$-weak attraction of $\mu$) has positive Lebesgue measure:


  \begin{equation} \label{equationbasinweakattraction} A_{\epsilon}(\mu):= \{ x \in S^1: \mbox{dist}(p \omega (x), \mu) < \epsilon\}.\end{equation}
  \end{proposition}


  For the sake of completeness, and although Proposition \ref{propositionPhysicalLike} can be easily obtained from the results in \cite{CE}, we give an independent proof   in Section 3 of this paper. Let us state now our main result:

\begin{theorem}
\label{maintheorem}
For any $C^1$-expanding map $f: S^1 \mapsto S^1$ there exist   SRB-like  measures and all of them are   equilibrium states for the potential $- log |f'|$.   \em

\vspace{.3cm}

The following assertions are immediate  consequences or  restatements  of Theorem \ref{maintheorem},  for  all the $C^1$ expanding maps that are non necessarily $C^{1+ \alpha}$:

\vspace{.3cm}

\noindent {\bf \ref{maintheorem}.1} \em The set $ES_f$ of equilibrium states for $-\log |f'|$ contains the weak$^*$-compact convex hull of the  never empty  set ${\mathcal O}_f$ of SRB-like  measures.

\vspace{.3cm}

\noindent \em {\bf \ref{maintheorem}.2} \em If   $ES_f$ contains a single measure $\mu$, then $\mu$ is ergodic and ${\mathcal O}_f = \{\mu\}$. Besides,   ${\mathcal O}_f = \{\mu\}$ if and only if $\mu$ is SRB and its basin $B(\mu)$ has full Lebesgue measure. \em

\vspace{.3cm}

\noindent  {\bf \ref{maintheorem}.3} \em  Any SRB-like measure $\mu $   (and in particular any SRB measure if it exists) satisfies Pesin's Entropy Formula \em  (\ref{equationPesinsFormula}).

\vspace{.3cm}

\noindent {\bf \ref{maintheorem}.4} \em There  exist $f$-invariant probability measures such that $m(A_{\epsilon}(\mu)) >0$ for all $\epsilon >0$, where $m$ denotes the Lebesgue measure and $A_{\epsilon}(\mu)$ denotes the  basin of $\epsilon$-weak attraction of $\mu$    defined by \em (\ref{equationbasinweakattraction}). \em All those measures satisfy Pesin's Entropy Formula. \em

\end{theorem}

We prove Theorem \ref{maintheorem} in Section 4.
It  is a stronger version of Theorem 6.1.8 of the book of Keller \cite{ke}, that states that observable measures belong to $ES_f$. In fact, the definition
of SRB-like measures in Section 2 of this paper  is  non trivially weaker than the definition of observable measures in \cite{ke}. While SRB-like measures do exist for any $f \in {\mathcal E}^{1}$, the stronger observable measures according to \cite{ke} may not exist. Nevertheless, some of the arguments that we use to prove    Theorem \ref{maintheorem}, are taken from  the proof of Theorem 6.1.8 in \cite{ke}. The   difference   resides in the proof of Lemma  \ref{lemma1teoremon}. We have to manage with sets of probabilities (neighborhoods of the SRB-like measures), instead of fixed probabilities (the observable measures according to \cite{ke}).

The statement $\ref{maintheorem}$.2 can be  equivalently reformulated, substituting the assumption   $\#ES_f = 1$ by the following condition (see Lemma 2.4 \cite{qiuhao}):


$$\lim_{t \rightarrow 0^+} \frac{1}{t} \sup_{\nu \in {\mathcal M}_f} \Big(h_{\nu} + \int (t\varphi - \log |f'|) \, d \nu \Big) = \int \varphi \, d \mu \ \ \forall \ \varphi \in C^0(S^1)  \ \ \forall \ \mu \in ES_f.$$


For any expansive map (in any finite-dimensional manifold) the above condition is $C^1$ generic (see Corollary 2.5 and Proposition 3.1 of \cite{qiuhao}). Thus, the statement \ref{maintheorem}.2 provides a new proof of a remarkable result in \cite{cq}: $C^1$-generically the expanding maps of the circle have a unique ergodic SRB measure whose basin covers Lebesgue-almost all the orbits.

 Let us state three   Corollaries of Theorem \ref{maintheorem}. We say that a probability measure is atomic if it is supported on a finite set.

\begin{corollary}.  \label{corolarioteoremon}

There is no atomic   SRB-like measure of a $C^1$ expanding map in $S^1$.
\end{corollary}


We prove this Corollary in the paragraph \ref{proofCorollary1}.


\begin{corollary} \label{corollarymu<<m}
Denote by $m$   the Lebesgue measure in $S^1$. For any $C^1$ expanding map $f$ in $S^1$  the following assertions are equivalent:

\begin{description}
\item
\em (a) \em There exists some SRB-like measure $\mu$ such that $m \ll \mu$

\item \em (b) \em There exists a unique SRB-like measure $\mu$,  it is equivalent to Lebesgue and ergodic.
\end{description}
\noindent Besides, if the conditions above hold, then $\mu$ is SRB and its basin  $B(\mu)$ has full Lebesgue measure.

\end{corollary}


We prove this Corollary in the paragraph \ref{proofCorollary1} at the end of this paper. This corollary has a similar version for natural measures, when they exist, instead of SRB-like measures (Theorem 2.4, Part (3) of \cite{jarvenpaatolonen}).
From the definition of SRB-like measure, it is immediate that if there exists some ergodic SRB-like measure $\mu$ such that $\mu \ll m$, then it is SRB. Nevertheless it may exist non ergodic invariant measures $\mu \ll m$ that are neither SRB nor SRB-like (see \cite{Q3}). In such a case  $\mu$ satisfies Pesin's Entropy Formula, as stated in the following lemma.  This shows that the SRB-like condition is sufficient but not necessary to a measure $\mu$ be an equilibrium state for the potential $-\log |f'|$.


\begin{corollary}
\label{corolario3}
  Let $f$ be a $C^1$ expanding map of $S^1$. Let $\mu$ be a non ergodic  $f$-invariant probability such that $\mu \ll m$, where $m$ is the Lebesgue measure. Then $\mu$  satisfies Pesin's Entropy Formula.

\end{corollary}


The proof of Corollary \ref{corolario3} is in the paragraph \ref{proofCorollary3}. This corollary has a similar formulation
for $C^1$-diffeomorphisms in any dimension with a dominated splitting (see \cite{suntian}).

\section{SRB-like measures.}

  We revisit the  definition and properties  of  the  SRB-like (weakly physical) measures. The content of this section is a reformulation   of a  part of      \cite{CE}.

 \begin{proposition} \label{PropositionObservable}
 There exists a  unique  minimal  nonempty and weak$^*$ compact set ${\mathcal O}_f \subset {\mathcal M}$ such that $p\omega(x) \subset {\mathcal O}_f$ for a full-Lebesgue set of initial states $x \in S^1$.
 \end{proposition}


 {\em Proof:} Consider the family ${ \Upsilon }$ of all the non empty and weak$^*$ compact sets ${\mathcal A} \subset {\mathcal M}$ such that $p\omega(x) \subset {\mathcal A}$ for a full Lebesgue set of initial states $x \in S^1$. The family ${\Upsilon}$ is not empty, since trivially ${\mathcal M} \in {\Upsilon}$. Define in ${\Upsilon}$ the partial order ${\mathcal A}_1 \leq {\mathcal A}_2$ if and only if ${\mathcal A}_1 \subset {\mathcal A}_2 $.
 We assert that each chain in $\Upsilon$ has a minimal element in $\Upsilon$. In fact, $\{{\mathcal A}_{\alpha}\}_{\alpha \in \aleph} \subset \Upsilon$ is a chain if it is a totally ordered subset of $\Upsilon$. Let us prove that ${\mathcal A}:= \bigcap_{\alpha \in \aleph} {\mathcal A}_{\alpha} $ belongs to $\Upsilon$. For each fixed $\alpha \in \aleph$,  and for each $\epsilon >0$   define
 $B_0 ( \alpha) := \{x \in M: p\omega(x) \subset {\mathcal A}_{\alpha} \}$, $B_{\epsilon}( {\mathcal A}) := \{x \in M: p \omega(x) \subset {\mathcal B}_{\epsilon}({\mathcal A}  )\}$, where ${\mathcal B}_{\epsilon}({\mathcal A}) := \{\nu \in {\mathcal M}: \mbox{dist}(\nu, {\mathcal A}) < \epsilon\}$. To conclude that ${\mathcal A} \in \Upsilon$, it is enough to prove that $m (B_{\epsilon}( {\mathcal A})) = 1$ for all $\epsilon >0$, where $m $ denotes the Lebesgue measure on $M$.
  For all $\epsilon >0$ there exists $\alpha \in \aleph$ such that ${\mathcal A}_{\alpha} \subset {\mathcal B}_{\epsilon}({\mathcal A})$. (If it did not exist  then, by the property of finite intersections of compact sets, and since $\{{\mathcal A}_{\alpha}\}_{\alpha \in \aleph}$ is totally ordered, we would deduce that the set $\bigcap_{\alpha \in \aleph} \big ({\mathcal A}_{\alpha} \setminus {\mathcal B}_{\epsilon}({\mathcal A})\big) $ would be nonempty, contained in ${\mathcal A}$, but disjoint with its open neighborhood ${\mathcal B}_{\epsilon}({\mathcal A})$.) We deduce that $B_0(\alpha) \subset B_{\epsilon}( {\mathcal A})$. Since $A_{\alpha} \in \Upsilon$, we have that  $m(B_0({\alpha}))= 1$ for all $\alpha \in \aleph$. Thus $m(B_{\epsilon}( {\mathcal A})) = 1$ for all $\epsilon >0$, and therefore ${\mathcal A} \in \Upsilon$. We have proved that each chain in $\Upsilon$ has a minimal element in $\Upsilon$. So, after Zorn Lemma there exist minimal elements in ${\Upsilon}$, namely, minimal non empty and weak$^*$ compact sets ${\mathcal O} \subset {\mathcal M}$ such that $p\omega(x) \subset {\mathcal O}$ for Lebesgue almost all $x \in S^1$. Finally, the minimal element ${\mathcal O} \subset {\Upsilon}$ is unique since the intersection of two of them is also in ${\Upsilon}. \ \ \ \ \Box$

\begin{definition}
 \em \label{definitionobservable}  {\bf (SRB-like probability measures.)}
A probability measure $\mu \in \mathcal M$ is \em SRB-like \em or \em weakly physical \em
if  $\mu \in {\mathcal O}_f$, where ${\mathcal O}_f$ is the set of Proposition \ref{PropositionObservable}.

\end{definition}


 It is immediate  that any SRB-like measure is $f$-invariant. In fact, the set of $f$-invariant Borel probabilities is non empty, weak$^*$-compact and contains $p\omega(x)$ for all $x \in S^1$.
It is also immediate that all the SRB measures (according with Definition \ref{definitionPhysical}), if they exist, are SRB-like measures. In fact, if  $\mu \not \in {\mathcal O}_f$, then  since $p \omega(x) \subset{\mathcal O}_f $ for Lebegue-almost all $x \in S^1$, the set $B(\mu)= \{ x \in S^1: p\omega(x)= \{\mu\}\}$  has zero Lebesgue-measure,  and thus $\mu $ is non SRB.

  \begin{prueba}
  \label{proofPropositionPhysicalLike}
  {\bf Proof of   Proposition \ref{propositionPhysicalLike}}

    \em As said in Section \ref{sectionStatements}, this Proposition    gives an individual (weakly) physical meaning to each of the SRB-like measures.

{\em Proof:} Let us denote with $m$ the Lebesgue measure. For any $\epsilon >0$ and any $\mu \in {\mathcal M}$ let us denote ${\mathcal B}_{\epsilon}(\mu)$ to the ball of center $\mu$ and radius $\epsilon$ in ${\mathcal M}$, defined with the metric ${\mbox{dist}}$.
If $\mu \in {\mathcal O}_f$ then      $m(A_{\epsilon}(\mu)) >0$ for all $\epsilon >0$, because if not, the compact set ${\mathcal K}:= {\mathcal O}_f\setminus{\mathcal B}_{\epsilon}(\mu)$ would be strictly contained in ${\mathcal O}_f$ and such that $p \omega(x) \subset {\mathcal K}$ for Lebesgue almost all $x \in S^1$. (Therefore ${\mathcal K}$ is not empty.)  This last contradicts the minimality condition of ${\mathcal O}_f$ in   Proposition \ref{PropositionObservable}. Thus any SRB-like measure satisfies the statement $m(A_{\epsilon}(\mu)) >0$ for all $\epsilon >0$.
On the other hand, if  a Borel probability measure $\mu$ satisfies the inequality $m(A_{\epsilon}(\mu)) >0$   for all $\epsilon >0$, and  since $p \omega(x) \subset {\mathcal O}_f$ for $m$ a.e. $x \in S^1$, we obtain  ${\mathcal B}_{\epsilon}(\mu) \bigcap {\mathcal O}_f \neq \emptyset$ for all $\epsilon >0$. Namely, $\mu$ is in the weak$^*$-closure of   ${\mathcal O}_f$. Since ${\mathcal O}_f$ is weak$^*$-compact (see Proposition \ref{PropositionObservable})  we conclude that $\mu \in {\mathcal O}_f$ as wanted. $\ \ \ \ \Box$

\end{prueba}

\section{Proof  of Theorem \ref{maintheorem}}

Denote by ${\mathcal M}_ f \subset {\mathcal M}$ the set of all the $f$-invariant Borel probability measures on $S^1$.
Let us recall the definition of the measure theoretic entropy.
For any Borel measurable finite partition ${\mathcal P}$ of $S^1$, and for any (non necessarily invariant) probability $\mu$ it is defined


$$H({\mathcal P}, \mu) := -\sum_{X_i \in {\mathcal P}} \mu(X_i) \log \mu(X_i) $$
If besides $\mu \in {\mathcal M}_f$ then $$  h ({{\mathcal P}, \mu}):= \lim_{q \rightarrow + \infty} \frac{H({\mathcal P}^q , \mu)}{q}$$
\noindent In the equality above ${\mathcal P}^q := \bigvee_{j= 0}^{q-1} f^{-j}({\mathcal P})$, where for any pair of  finite partitions ${\mathcal P}$ and ${\mathcal Q} $ it is defined  ${\mathcal P}\bigvee {\mathcal Q}:= \{X  \bigcap Y  \neq \emptyset: \ X  \in {\mathcal P}, \ Y \in {\mathcal Q}    \}$. It is a well established result that the limit defining $h ({{\mathcal P}, \mu})$ exists.
 Finally, the measure theoretic entropy $h_{\mu}$ of an $f$-invariant measure $\mu$ is defined by
$h _{\mu} := \sup _{\mathcal P} h (\mu, {\mathcal P}),$
where the sup is taken on all the Borel measurable finite partitions ${\mathcal P}$ of the space.

We define the diameter $\mbox{diam}{\mathcal P}$ of a finite partition ${\mathcal P}$ as the minimum diameter of its pieces. A well known result (see Proposition 2.5 of \cite{Bowen0}) states that if $f$ is expansive (in particular if $f \in {\mathcal E}^1$), and if ${\mathcal P}$ is a partition with diameter smaller than the expansivity constant,  then $h_{\mu} = h(\mu, {\mathcal P})$ for all $\mu \in {\mathcal M}_f$. Applying this result, in the sequel we will consider only finite partitions with diameter smaller than the expansivity constant $\alpha$ of $f \in {\mathcal E}^1$. So, we will compute the measure theoretic entropy by


\begin{equation}
\label{equationEntropy}
h_{\mu} = \lim_{q \rightarrow + \infty} \frac{H({\mathcal P}^q , \mu)}{q} \ \ \mbox{ if } \ \  \mbox{diam}({\mathcal P}) < \alpha.
\end{equation}


\noindent For any (non necessarily $f$-invariant) Borel probability $\nu$,   denote $f^* \nu$ to the probability defined by $f^* \nu (B) := \nu(f^{-1}(B))$ for any Borel-measurable set $B$. For a given  finite partition ${\mathcal P}$ denote $\partial {\mathcal P} := \bigcup _{X   \in {\mathcal P}} \partial X$, where $\partial X$ denotes the topological boundary of the piece $X$. The only step along the proof of Theorem \ref{maintheorem} (which is one of the key-points of this proof), for which we use that the space has dimension one, resides in the application of the following lemma, in particular in its statements (ii) and (iii). This lemma is essentially a  restatement of a part of Misiurewicz's proof of the Variational Principle:
 \begin{lemma}
 \label{lemmaEntropy} Let $f $ be a $C^1$ expanding map  on $S^1$. Let $\alpha >0$ be an expansivity constant. For all $0 < \delta \leq \alpha$ there exists a finite partition ${\mathcal P}$ of $S_1$ such that:

\noindent \em (i)    $\mbox{diam}({\mathcal P}) < \delta \leq \alpha  $, \em

\noindent \em (ii) \em $\mu (\partial {\mathcal P} ) = 0$ for all $\mu \in {\mathcal M}_f$,

\noindent \em (iii) \em For any sequence of non necessarily invariant probabilities $\nu_n$,
 for any $\mu \in {\mathcal M}_f$  equal to the weak$^*$ limit of a convergent subsequence $\{\mu_{n_i}\}_{i \geq 1}$ of $\mu_n:= \frac{1}{n} \sum_ {j= 0}^{n-1} (f^j)^* \nu_n$,
and for any $\epsilon >0$, there exists $i_0$ such that
$$ \frac{1}{n_i} H({\mathcal P}^{n_i}, \nu_{n_i}) \leq h_{\mu} + \epsilon \ \ \ \forall \ i \geq i_0.$$
\em

 \end{lemma}


{\em Proof: } Take any finite covering ${\mathcal U}$ of $S^1$ with open intervals with length smaller than $\delta$.  Denote $\partial{\mathcal U} := \bigcup_{X \in {\mathcal U}} \partial X  $. It is a finite set. Therefore $\mu (\partial {\mathcal U}) = 0$ for all $\mu \in {\mathcal M}_f$  if and only if $\partial {\mathcal U}$ does not contain periodic points of $f$. Since $f \in {\mathcal E}^1$, the set of periodic points is countable. Then, changing if necessary the open intervals $X \in {\mathcal U}$ to slightly smaller ones such that they still cover $S^1$ and their boundary points are non periodic, we get a new covering  ${\mathcal U}' = \{  Y_i\}_{1 \leq i \leq p}$   such that $\mu (\partial {\mathcal U}') = 0$ for all $\mu \in {\mathcal M}_f$. Therefore, the partition ${\mathcal P} = \{X_i\}_{1 \leq i \leq p}$ defined by $X_1 := Y_1 \in {\mathcal U}'$, \  $X_{i+1} :=  Y_{i+1} \setminus (\cup_{j= 1}^i X_i)$, satisfies the assertions (i) and (ii).
Let us prove that   (i) and (ii) imply (iii).
Fix the integer numbers $q \geq 1$,  and $n \geq q$. Write $n= N q + j$ where $N, j$ are integer numbers such that $ 0 \leq j \leq q-1$ Fix a (non necessarily invariant) probability $\nu$. From the properties of the entropy function $H$ of $\nu$ with respect to the partition ${\mathcal P}$, we obtain


$$H({\mathcal P}^n, \nu) = H({\mathcal P}^{ Nq+ j }, \nu)   \leq H( \vee_{i= 0}^{j-1}  f^{-i}{\mathcal P}, \nu) + H({\vee_{i= 1}^{N}f^{-iq}\mathcal P}^{q }, \nu)  \leq $$


$$\sum_{i= 0}^{j-1} H(f^{-i} {\mathcal P}, \nu) + \sum_{i= 1}^N H(f^{-iq}{\mathcal P}^{q },  \nu) = \sum_{i= 0}^{j-1} H({\mathcal P} , (f^{i})^* \nu) + \sum_{i= 1}^N H({\mathcal P}^{q }, (f^{iq})^*\nu)   $$


$$\Rightarrow \ H({\mathcal P}^n, \nu) \leq q \log p + \sum_{i= 1}^N H({\mathcal P}^{q }, (f^{iq})^*\nu) \ \ \forall \ q \geq 1, \ \ n \geq q.  $$


 \noindent  To obtain the inequality above
 recall that $H({\mathcal P}, \nu) \leq \log p \ \ \forall \ \nu \in {\mathcal M}$, where $p   $ is the number of pieces of the partition ${\mathcal P}$. The inequality above holds also for $f^{-l}{\mathcal P}$ instead of ${\mathcal P}$, for any $l \geq 0$, because it holds for any partition with exactly $p$ pieces.  Thus:


$$\ H(f^{-l}{\mathcal P}^n, \nu) \leq q \log p + \sum_{i= 1}^N H(f^{-l}{\mathcal P}^{q }, (f^{iq})^*\nu) =$$


$$ q \log p + \sum_{i= 1}^N H({\mathcal P}^{q }, (f^{iq + l})^*\nu). $$


\noindent Adding the above inequalities   for $0 \leq l \leq q-1$, we obtain:


$$\sum_{l= 0}^{q-1} H(f^{-l}{\mathcal P}^n, \nu) \leq  q^2 \log p + \sum_{l= 0}^{q-1} \sum_{i= 1}^N H({\mathcal P}^{q }, (f^{iq + l})^*\nu) $$


\begin{equation} \label{equation0}\Rightarrow \ \sum_{l= 0}^{q-1} H(f^{-l}{\mathcal P}^n, \nu) \leq q^2 \log p + \sum_{l= 0}^{Nq+ q -1} H({\mathcal P}^{q }, (f^{l })^*\nu). \end{equation}


\noindent On the other hand, for all $0 \leq l \leq q-1$


$$H({\mathcal P}^n, \nu) \leq H({\mathcal P}^{n+ l}, \nu) \leq \Big (\sum_{i= 0}^{l-1} H(f^{-i}{\mathcal P}, \nu) \Big) + H(f^{-l} {\mathcal P}^n, \nu)     $$


$$\leq q\log p + H(f^{-l} {\mathcal P}^n, \nu).$$


\noindent Therefore, adding the above inequalities   for $0 \leq j \leq q-1$ and joining with the inequality (\ref{equation0}), we obtain:


  $$q H({\mathcal P}^n, \nu)  \leq 2 q^2 \log p + \sum_{l= 0}^{Nq+ q -1} H({\mathcal P}^{q }, (f^{l })^*\nu).$$
Recall that $n = Nq + j$ with $0 \leq j \leq q-1$. So $Nq + q \leq n + q$ and then


$$q H({\mathcal P}^n, \nu)  \leq 2 q^2 \log p \ + \sum_{l= 0}^{n -1} H({\mathcal P}^{q }, (f^{l })^*\nu) \ + \sum_{l=n}^{Nq + q-1} H({\mathcal P}^q, (f^{l })^*\nu) $$


$$\Rightarrow \ \ q H({\mathcal P}^n, \nu) \leq  3 q^2 \log p + \sum_{l= 0}^{n -1} H({\mathcal P}^{q }, (f^{l })^*\nu) .$$


\noindent  In the last inequality we have used that the number of nonempty pieces of ${\mathcal P}^q$ is at most $p^{q}$. Now we put $\nu= \nu_n$ and divide by $n$. Recall  that the convex combination of the function $H$ for a finite set of probability measures is not larger than the function $H$ for the convex combination of the measures. We deduce:


  $$ \frac{q\, H({\mathcal P}^n, \nu_n)}{n} \leq  \frac{3 q^2 \log p}{n} + \frac{1}{n} \, \sum_{l= 0}^{n -1} H({\mathcal P}^{q }, (f^{l })^*\nu_n)  $$


  $$\Rightarrow \ \   \frac{q\, H({\mathcal P}^n, \nu_n)}{n}  \leq \frac{3 q^2 \log p}{n} + H\Big({\mathcal P}^{q }, \mu_n \Big) .$$


   \noindent For any fixed $\epsilon >0$ (and the natural number $q \geq 1$ still fixed), take    $n \geq n(q):= \max\{q, 9\, q  \log p/ \epsilon \}$ in the inequality above. We deduce:


 $$\frac{q}{n} H({\mathcal P}^n, \nu_n) \leq  \frac{q \epsilon}{3} + H({\mathcal P}^q , \mu_n) \ \ \ \ \forall \ n \geq n(q)  \ \ \ \forall \ q \geq 1. $$


\begin{equation}
\label{equationSuper}
\Rightarrow \ \ \frac{1}{n} H({\mathcal P}^n, \nu_n) \leq    \frac{\epsilon}{3} + \frac{H({\mathcal P}^q, \mu_n)}{q} \ \ \ \ \ \ \forall \   n \geq n(q) \ \ \ \ \forall \ q \geq 1.\end{equation}


\noindent The inequality above holds for  for any fixed $q \geq 1$ and for any $n$ large enough, depending on $q$.

By hypothesis, $\mu$ is   $f$-invariant equal to the weak$^*$-limit of a convergent subsequence of $\mu_n$. After Equality (\ref{equationEntropy}) there exists $q \geq 1$ such that


\begin{equation}
\label{equationSuper2}
\frac {H({\mathcal P}^q, \mu)}{q} \leq  h_{\mu} + \frac{ \epsilon}{3}.\end{equation}


\noindent Fix such a value of $q$. Since $\mu(\partial({\mathcal P})) = 0$ for all $\mu \in {\mathcal M}_f$:
$$\lim_{i \rightarrow + \infty} H({\mathcal P}^q, \mu_{n_i}) = H({\mathcal P}^q, \mu) \mbox{ if } \lim _{i \rightarrow + \infty} \!\!(\mbox{weak}^*) \ \mu_{n_i} = \mu,  \ \ \mbox{ because } \mu(\partial {\mathcal P}^q) = 0. $$
Therefore, there exists $i_0$ such that for all $i \geq i_0$:
$$n_i \geq n(q) \ \mbox{ and } \ \ \frac{H({\mathcal P}^q, \mu_{n_i}) }{q} \leq \frac{H({\mathcal P}^q, \mu)}{q} + \frac{\epsilon}{3}.$$
Joining the last assertion with Inequalities (\ref{equationSuper}) and (\ref{equationSuper2})  we deduce (iii), as wanted. $\ \Box$


\begin{prueba}
{\bf Notation:}  \em \label{notacionKsubr}  For any $f \in {\mathcal E}^1$ denote $\psi := - \log |f'| < 0$, and
 for all $r \geq 0$ construct


 \begin{equation}\label{equationKr}{\mathcal K}_r := \{\nu \in {\mathcal M}_f: \;
\int \psi \, d \nu + h_{\nu} \geq - r\}.\end{equation}


\noindent The notation above is  taken from   the book \cite{ke}. The set ${\mathcal K}_r$ is non empty, weak$^*$ compact  and convex. In fact, join
the proof of
Theorem
4.2.3  of the book in \cite{ke}     with  Theorem 4.2.4 and Remark  6.1.10
of the same book. Joining the assertions (\ref{equationRuelleInequality}), (\ref{equationRLYequality}) and (\ref{equationKr})  we deduce that $ES_f= {\mathcal K}_0$ is weak$^*$ compact and convex.
\end{prueba}


For any integer $n \geq 1$  and for all $x \in S^1$ recall the definition of the empirical probability $\sigma_n(x)$ in Equality (\ref{eq1}), and the definition of the p-limit set $p\omega(x)$ in the set  $  {\mathcal M}$ of Borel probabilities, according to   Equality (\ref{equationpomegax}).
In  ${\mathcal M}$   fix the following weak$^*$ metric:


  \begin{equation}
  \label{equationweak^*distance}
  \mbox{dist}(\mu, \nu) := \sum_{i= 0}^{+ \infty} \frac{1}{2^i} \; \left |\int \phi_i \, d \mu - \int \phi_i \, d \nu\right |, \end{equation}


 \noindent where $\phi_{\,0}:= \psi= - \log |f'|$ and $\{\phi_i\}_{i \geq 1}$ is a countable family of continuous functions   that is dense in the space $C^0(S^1, [0,1])$. Trivially  with this distance, for any $\mu_0 \in {\mathcal M}$ and any $\epsilon >0$ the ball ${\mathcal B} := \{\nu \in {\mathcal M}: \ \mbox{dist} (\mu_0, \nu) < \epsilon\}$   is convex.

   \begin{lemma}
\label{lemma1teoremon}  Let $f$ be a $C^1$ expanding map on
  $S^1$. Let $m$ be the Lebesgue measure on $S^1$. Fix $r >0$ and let ${\mathcal K}_r  $ be defined by Equality \em (\ref{equationKr}). \em Consider the weak$^*$ distance defined in \em (\ref{equationweak^*distance}). \em Then, for all $0 < \epsilon < r/2$ there exists $n_0 \geq 1$  such that  \em \begin{equation} \label{equationlemma1teoremon}
       m (\{x \in S^1: \mbox{dist} ( \sigma_{n}(x)  , {\mathcal K} _r )  \geq \epsilon\} ) \leq   e^{ n(\epsilon  - r)} <  e^{- {nr}/{2}}  \ \ \ \ \ \forall \ n \geq n_0.\end{equation}

\end{lemma}

{\em Proof: }
   Fix $0 <\epsilon< r/2$. Observe that the set $ \{\mu \in {\mathcal M} : \mbox{dist}(\mu, {\mathcal K}_r) \geq  \epsilon\}$ is weak$^*$ compact, so it has a finite   covering $\{{\mathcal B}_i\}_{1 \leq i \leq k} $, for a minimal cardinal $k \geq 1$, with open balls ${\mathcal B}_i \subset {\mathcal M}$ of radius $  \epsilon/3 $. For any fixed $n \geq 1$ denote


     $$C_{n,i} := \{x \in S^1: \ \sigma_{n}(x) \in {\mathcal B}_i\}, \ \ C_n := \bigcup_{i= 1}^k C_{n,i}. \ \ \mbox{ Then: } $$


 $$\{x \in S^1: \ \ \mbox{dist}(\sigma_{n}(x), {\mathcal K}_r) \geq \epsilon \} \subset C_n . $$ Therefore, to prove the lemma it is enough to find $n_0$ such that $m(C_n) \leq e ^{n(\epsilon - r)}$ for all $ n \geq n_0$. Fix $1 \leq i \leq k$. We claim:
\begin{equation}
\label{equationToBeProved}
  \exists \ n_i\ \mbox{ such that } \  m(C_{n,i}) \leq e ^{n(  - r + (\epsilon/2))} \ \ \forall \ \ n \geq n_i \ \ \ \mbox{ (to be proved)} \end{equation}

First, let us see that it is enough to prove Assertion (\ref{equationToBeProved}) to end the proof of the lemma. In fact,  if Assertion (\ref{equationToBeProved}) holds,   define:

\noindent $n_0:= \max\{ 2(\log k )/\epsilon, \max _{i= 1}^k n_i\}. $  Then, we deduce the following inequalities for all $n \geq n_0$, as wanted:


$$m(C_n) \leq \sum_{i= 1}^k m (C_{n,i}) \leq k \ e^{n(-r + \epsilon/2))} = e^{n(-r + (\epsilon/2) + (\log k /n))} \leq e^{n (\epsilon - r)}. $$

Second and last, let us prove   Assertion (\ref{equationToBeProved}).
 Consider an expansivity constant $\alpha >0$ of $f$. Take $\epsilon /6$ and for such value, fix  a continuity modulus $0 <\delta < \alpha$ of the function $\psi = - \log |f'|$. Namely  $|\psi(x) - \psi (y)| < \epsilon/6$ if $\dist (x,y) < \delta$. Take a finite partition ${\mathcal P} = \{X_i\}_{1 \leq i \leq p}$ of $S^1$  with diameter smaller than $ \delta  $ and satisfying also the conditions (ii) and (iii) of Lemma \ref{lemmaEntropy}. The   map $f$  is conjugated to a  linear expanding map. Therefore, if the diameter of the partition ${\mathcal P}$ is chosen small enough, the restricted map  $f^n|_{X} : X \mapsto f^n(X)$ is a diffeomorphism for all  $X \in {\mathcal P}^n$ and for all $n \geq 1$. Thus, recalling that $\psi= - \log |f'|$, we deduce the following equality for all $X \in {\mathcal P}^n$:
 $$m(X \cap C_{n, i}) = \int _{f^n(X \cap C_{n,i})} |(f^{-n})'| \, d m = \int _{f^n(X \cap C_{n,i})} e ^{\sum_{j= 0}^{n-1} \displaystyle{\psi \circ f^j }} \, d m  $$
 $$\Rightarrow \ \ \ m(C_{n, i})  = \sum_{X \in {\mathcal P}^n} \int _{f^n(X \cap C_{n,i})} e ^{\sum_{j= 0}^{n-1}\displaystyle{ \psi \circ f^j } }\, d m.$$
 Either $C_{n,i} = \emptyset$, and Assertion (\ref{equationToBeProved}) becomes trivially proved, or the finite family  of pieces $  \{X \in {\mathcal P}^n : \ X \cap C_{n,i} \neq \emptyset \} = \{X_1, \ldots, X_N\} $ has $N = N(n,i)$ pieces for some $N \geq 1$. In this latter case, take a unique point $y_k \in X_k  \cap C_{n,i}$ for each $k= 1, \ldots, N$. Denote by
 $Y(n,i) = \{y_1, \ldots, y_{N}\}   $  the collection of such points. Due to the construction of $\delta >0$, and since the partition ${\mathcal P}$ has diameter smaller than  $\delta$, we deduce:

  $$\sum_{j= 0}^{n-1} \psi(f^j(y)) \leq \sum_{j= 0}^{n-1} ( \psi (f^j(y_k)) +   \epsilon/6)   \ \ \forall \ y, y_k \in X_k, \ \forall \  k= 1, \ldots, N.$$


 \noindent Therefore
$m(C_{n, i})  \leq e^{n \epsilon/6} \sum_{k= 1}^N e ^{\sum_{j= 0}^{n-1} \displaystyle{\psi (f^j(y_k))}} \, m(f^n(X_k \cap C_{n,i}))   $, and thus:


  $$m(C_{n, i}) \leq e^{n \epsilon/6} \sum_{k= 1}^N e ^{\sum_{j= 0}^{n-1} \displaystyle{\psi (f^j(y_k))}}.$$
\noindent Define


$$\ \ \ \ \ \ \ \ \ \ L:= \sum_{k= 1}^N e ^{\sum_{j= 0}^{n-1} \displaystyle{\psi (f^j(y_k))}}, \ \ \ \ \ \lambda_k : = \frac{1}{L} \, e ^{\sum_{j= 0}^{n-1} \displaystyle{\psi (f^j(y_k))}} >0. $$  Then, $ \sum_{k= 1}^N \lambda_k = 1$ and


 $$m(C_{n, i}) \leq e^{(n \epsilon/6) + \log L} , \ \ \ \log L = \left (\sum_{k= 1}^N \lambda_k \sum_{j= 0}^{n-1} \displaystyle{\psi (f^j(y_k))} \right ) - \left (\sum_{k= 1}^{N} \lambda_k \log \lambda_k \right ).$$
 (To check the last equality  substitute $\lambda_k$ by the quotient  which defines it.)

  \noindent Define the probability measures


  $$\nu_{n} := \sum_{k= 1}^N \lambda_k \delta_{y_k}, \ \ \ \ \mu_{n } := \frac{1}{n} \sum _{j= 0}^{n-1} (f^j)^* (\nu_n) = \sum_{k=1}^N \lambda _k \frac{1}{n} \sum_{j= 0}^{n-1} \delta_{f^j(y_k)} = \sum_{k= 1}^N \lambda_k\sigma_{n}(y_k).  $$

  \noindent It is standard to check that


$$\sum_{k= 1}^N \lambda_k \sum_{j= 0}^{n-1} \displaystyle{\psi (f^j(y_k))} = n\int \psi \, d \mu_n, \ \ \ \ \sum_{k= 1}^{N} \lambda_k \log \lambda_k  = H({\mathcal P}^n, \nu_n), \ \ \mbox{ and then} $$
  $$m(C_{n, i}) \leq \mbox{exp}\Big({\frac{n \epsilon}{6} + \log L}\Big) = \mbox{exp} \Big ( n \Big ( \frac{\epsilon}{6} + \int \psi \, d \mu_n + \frac{H({\mathcal P}^n, \nu_n)}{n} \Big ) \Big ) $$
  Take a subsequence $n_l \rightarrow + \infty$ such that

  $\bullet$ $\lim_{l \rightarrow + \infty} \frac{1}{n_l} \log m(C_{n_l, i}) = \limsup _{n \rightarrow + \infty} \frac{1}{n} \log m (C_{n, i})$ and

   $\bullet$ the sequence   $\{\mu_{n_l}\}_{l \geq 1}$ is weak$^*$-convergent.

   \noindent Denote   $\mu = \lim_ {l \rightarrow + \infty} \mu_{n_l}$.
   After Assertion (iii) of Lemma \ref{lemmaEntropy} and   the definition of the weak$^*$ topology, there exists $n_i \geq 1$ such that
  \begin{equation} \label{equationPrincipal} m(C_{n,i}) \leq \mbox{exp} \Big ( n \Big ( \frac{\epsilon}{2} + \int \psi \, d \mu  + h _{\mu} \Big ) \Big ) \ \ \forall \ n \geq n_i. \end{equation}
  By construction  $y_k \in C_{n,i}$ for all $k= 1, \ldots, N$. Thus $\sigma_{n}(y_k) \in {\mathcal B}_i$. Since the ball ${\mathcal B}_i$ is convex and $\mu_n$ is a convex combination of the measures $\sigma_{n}(y_k)$ (recall that $ \sum_{k= 1}^N \lambda_k= 1$), we deduce that $\mu_n \in {\mathcal B}_ i $.  Therefore, the weak$^*$ limit $\mu$ of any convergent subsequence of $\{\mu_n\}_{n \geq 1}$ belongs to ${\overline{\mathcal B}_i}$. Since the ball ${\mathcal B}_i$ has radius $\epsilon/3$ and intersects $\{\mu \in {\mathcal M} : \ \ \mbox{ dist} (\mu, {\mathcal K}_r) \geq \epsilon\}$, we have $\mu \in \overline {{\mathcal B}_i} \subset {\mathcal M}\setminus {\mathcal K}_r $. Therefore, by the definition of the set ${\mathcal K}_r$, we have:
  $h_{\mu} + \int \psi \, d \mu  < -r$.
  Substituting this last inequality in (\ref{equationPrincipal}) we conclude (\ref{equationToBeProved}) ending the proof. $ \ \Box$

\begin{prueba}

\label{pruebadirectoTeoremon}
{\bf End of the proof of Theorem \ref{maintheorem}} \em

For any $r >0$ consider the   compact set
${\mathcal K}_r \subset {\mathcal M}$ defined by Equality (\ref{equationKr}).
Since $\{{\mathcal K}_r\}_r$ is decreasing with $r $:


$${\mathcal K}_0  = \bigcap_{r >0}{\mathcal K}_r $$
From Equalities (\ref{equationRuelleInequality}) and  (\ref{equationRLYequality})  and from the definition of ${\mathcal K}_0$ in Equality (\ref{equationKr}), we have


$${\mathcal K}_0 = ES_f.$$ So, to prove Theorem \ref{maintheorem} me must prove that the set ${\mathcal O}_f$ of SRB-like measures satisfy: ${\mathcal O}_f \subset {\mathcal K}_r$ for all $r >0$. Since ${\mathcal K} _r$ is weak$^*$ compact, we have   $${\mathcal K}_r = \bigcap_{  \epsilon >0} {\mathcal B}(r, \epsilon), \ \ \ \  \mbox{ where } {\mathcal B}(r, \epsilon):= \{\mu \in {\mathcal M}: \mbox{dist}(\mu, {\mathcal K}_r) \leq \epsilon\}, $$ with the weak$^*$ distance defined in (\ref{equationweak^*distance}). Therefore, it is enough to prove that ${\mathcal O}_f \subset {\mathcal B}(r, \epsilon)$ for all $0 <\epsilon < r/2$ and for all $r >0$.
After   Proposition \ref{PropositionObservable}, and since ${\mathcal B}(r, \epsilon)$ is weak$^*$ compact, it is enough to prove that the following set $B (r, \epsilon)$ (called basin of attraction of ${\mathcal B}(r, \epsilon)$) has full Lebesgue measure: $${B(r, \epsilon)}
:= \{x \in S_1: \   pw(x)  \subset {\mathcal B}(r, \epsilon)  \}.$$

 From Lemma \ref{lemma1teoremon}, there exists $n_0 $
 such that, for any $n>n_0$: $$m \{x:\sigma_{n}(x) \not \in {\mathcal B}(r, \epsilon) \}\leq
 e^{ n(\epsilon -  r) } \leq e^{-nr/2},$$ where $m$ denotes the Lebesgue measure. This implies that $$\sum_{n=1}^\infty m (x: \sigma_{n}(x)
 \not \in   {\mathcal B}(r, \epsilon)) < +\infty.$$  After   Borel-Cantelli Lemma it follows that
 $$m\left(\bigcap_{n_0=1}^\infty \bigcup_{n=n_0}^\infty \{x:\sigma_{n}(x) \not \in   {\mathcal B}(r, \epsilon) \}\right)=0.$$ In other words for $m$-a.e. $x \in S^1$,  there
 exists $n_0 \geq 1 $ such that $\sigma_{n}(x)\in {\mathcal B}(r, \epsilon) $ for all $n\geq n_0$. Hence,
 $pw(x) \subset {\mathcal B}(r, \epsilon) $ for
 $m$-almost all the points $x \in S^1$, as wanted.
$\Box$

\end{prueba}

\section{ Proofs of the Corollaries}

\begin{prueba} {\bf Proof of Corollary  \ref{corolarioteoremon} }
\label{proofCorollary1}
\em

   If $\mu $ is an atomic   invariant  measure for
an expanding map $f$, then $h_\mu(f)=0$. Since $ \psi = -\log f'   <0$,
we have $h_{\mu}(f) + \int \psi \, d \mu < 0$. Therefore $\mu$ does not satisfy Pesin's formula (\ref{equationRLYequality}).   After
 Theorem \ref{maintheorem}.3 the measure $\mu$ is not SRB-like. $\ \Box$
\end{prueba}

To prove Corollaries \ref{corollarymu<<m} and \ref{corolario3} we will use the following definition:

\begin{definition} \em \label{definitionErgodicSet}
For any $f$-invariant probability measure $\mu$   \em the weak$^*$-closure ${\mathcal K}(\mu)$ of the ergodic components of $\mu$, \em is the minimal nonempty and weak$^*$-compact   set of probabilities such that


\begin{equation} \label{equation99} \mu_x := \lim_{n \rightarrow + \infty} \frac{1}{n} \sum_{j= 0}^{n-1} \delta_{f^j(x)}  \in {\mathcal K}(\mu)\ \ \mu\mbox{-a.e. } \ x \in S^1.\end{equation}

\noindent After Birkhoff's Ergodic Theorem, the     above limit exists (in the weak$^*$ topology) $\mu$-a.e.   $x \in S^1$. Applying Zorn Lemma (as in the proof of Proposition \ref{PropositionObservable}, putting $\mu$ in the role of the Lebesgue measure $m$),  we deduce that the minimal compact set  ${\mathcal K}(\mu)$ satisfying (\ref{equation99}) exists and   is unique. We call ${\mathcal K}(\mu)$ the weak$^*$ closure of \em the ergodic components \em of $\mu$, because for $\mu$-a.e. $x \in X^1$ the limit $\mu_x$ in Equality (\ref{equation99}) is an ergodic component of $\mu$  (see for instance Theorem 4.1.12 of \cite{katokhasselblatt}.)
\end{definition}

\begin{lemma}
\label{lemmaMunotinErgodicSet}
For any $f$-invariant measure $\mu$, consider the weak$^*$-closure ${\mathcal K}(\mu)$ of its ergodic components, as defined in \em \ref{definitionErgodicSet}. \em Then, $\mu$ is ergodic if and only if $\mu \in {\mathcal K}(\mu)$, and if this latter inclusion occurs, then ${\mathcal K}(\mu) = \{\mu\}$. Thus, $\mu$ is non ergodic if and only if
$\mbox{dist} (\mu, {\mathcal K}(\mu)) >0$.
\end{lemma}

{\em Proof:} After Definition \ref{definitionErgodicSet}, and the definition of ergodicity, we have ${\mathcal K} = \{\mu\}$ if and only if $\mu$ is ergodic.  Now, it is enough to prove that if $\mu \in {\mathcal K}(\mu)$  then $\mu$ is ergodic. Consider the weak$^*$ distance defined by Equality (\ref{equationweak^*distance}). For any $\epsilon >0$ consider the ball ${\mathcal B}_{\epsilon} = \{\nu \in {\mathcal M}: \mbox{ dist}(\nu, \mu) < \epsilon\}$ and the set \begin{equation} \label{equationA} A_{\epsilon} = \{x \in S^1:  p \omega(x) \subset {\mathcal B}_{\epsilon} \}.\end{equation}  We claim that $\mu(A_\epsilon) >0$ for all $\epsilon >0$. In fact, arguing by contradiction if $\mu(A_{\epsilon}) = 0$, and since $p \omega (x) $ is a single measure for $\mu$-almost all the points $x \in S^1$, then $\mu( \{x \in S^1: p\omega(x) \subset {\mathcal K}(\mu)   \setminus {\mathcal B}_{\epsilon}\}) = 1$. This contradicts the minimality of   ${\mathcal K}(\mu)$   in Definition \ref{definitionErgodicSet}.
Consider the sequence of continuous functions $\phi_i$ in  Equality (\ref{equationweak^*distance}) which defines the weak$^*$ metric $\mbox{dist}$. Applying the Ergodic Decomposition Theorem (see for instance Theorem 4.1.12 of \cite{katokhasselblatt}):
   $$\int_{A_{\epsilon}} \phi_i d\mu = \int  \, d \mu  \int _{ A_{\epsilon}}  \phi_i \,  d \mu_x,$$
   where $\mu_x $ is an ergodic component of $\mu$. Since $A_{\epsilon}$ satisfies Equality (\ref{equationA}), and  $p\omega(x) = \{\mu_x\}$ for $\mu$-a.e. $x \in S^1$, we deduce $  \mu_x   \in {\mathcal B} _{\epsilon}$ for $\mu$-a.e. $x \in A_{\epsilon}$. Therefore
   $$\sum_{i= 0}^{+ \infty}    \frac{1}{2^i} \Big|\int \phi_i \, d \mu -  \int_{A_{\epsilon}} \phi_i d\mu\Big|     =  \sum_{i= 0}^{+ \infty}    \frac{1}{2^i} \Big | \int \phi_i \, d \mu -\int  \, d \mu  \int _{  A_{\epsilon}}  \phi_i  d \mu_x \Big |   $$ $$  \leq \sum_{i= 0}^{+ \infty}    \frac{1}{2^i} \int   \Big | \int  \phi_i \, d \mu - \int _{ A_{\epsilon}}  \phi_i  d \mu_x \Big | \, d \mu \leq 2\epsilon. $$
     The bounded linear operator $\varphi \in C^0(S^1, \mathbb{R}) \mapsto \int _{A_{\epsilon}} \varphi \, d \mu $ (via   Riesz Representation Theorem)  is the integral operator with respect to the finite measure $\mu_{\epsilon}$, defined by $$\mu_{\epsilon} (B) := \mu(A_\epsilon \cap B)$$ for all the Borel sets $B \subset S^1$.  The above inequality is translated as $$\mbox{dist}(\mu_{\epsilon}, \mu) \leq 2 \epsilon$$ in the space of all the finite Borel-measures $\nu $ such that $\nu (S^1) \leq 1$. Thus, $\lim_{\epsilon\rightarrow 0^+}\mu_{\epsilon} = \mu$ in such a space endowed with the weak$^*$-topology. We deduce that $\lim_{\epsilon \rightarrow 0^+} \int \varphi \, d \mu_{\epsilon} = \int \varphi \, d \mu$ for any $\varphi \in C^0(S^1, \mathbb{R}) $. In particular for the constant real function  $\varphi= 1$, we obtain that
      $\lim_{\epsilon \rightarrow 0^+} \mu(A_{\epsilon}) = 1$.
     Consider the basin $B(\mu)$ of attraction of $\mu$ defined by Equality (\ref{equationBasinOfAttraction}). By construction, the sets $A_{\epsilon}$  decrease to $B (\mu)$ when $\epsilon>0$ decrease going to zero. Thus,
     $B(\mu) = \bigcap_{\epsilon >0} A_{\epsilon}$ and $ \mu(B(\mu)) = \lim _{\epsilon \rightarrow 0^+} \mu (A_{\epsilon})  = 1$.
     Taking into account the definition of the basin $B(\mu)$  in Equality (\ref{equationBasinOfAttraction}), we conclude that $\lim_{n \rightarrow + \infty} (1/n) \sum_{j= 0}^{n-1} \delta_{f^j(x)} = \mu$ for $\mu$-a.e. $x \in S^1$. Thus  $\mu$ is ergodic, as wanted. $\ \ \Box$

\begin{prueba} {\bf Proof of Corollary   \ref{corollarymu<<m}}   \em
\label{proofCorollary2}

  Trivially, (b) implies (a) and it also implies that $\mu$ is SRB and that its basin of attraction has full-Lebesgue measure (recall Definitions  \ref{definitionPhysical} and \ref{definitionobservable}). So, it is only left to prove that
    (a) implies (b).

    Assume (a). Since $\mu$ is SRB-like, it is $f$-invariant. Using that $m \ll \mu$ and applying Birkhoff Theorem and Definition \ref{definitionErgodicSet}, we obtain that $p \omega(x) = \{\mu_x\}$  for $m$ a.e.   $x \in S^1$, where $\mu_x \in {\mathcal K}(\mu)$.  Applying Proposition \ref{propositionPhysicalLike} to the SRB-like measure $\mu$, for all $\epsilon>0$ there exists a $m$-positive set $A_{\epsilon}(\mu)$ such that $\mbox{dist}(p\omega(x), \mu) < \epsilon$. We deduce that $\mbox{dist}\big (\mu, {\mathcal K}(\mu) \big ) < \epsilon$ for all $\epsilon >0$. Thus, $\mu \in {\mathcal K}(\mu)$.  As proved in Lemma \ref{lemmaMunotinErgodicSet}, if $\mu$ were non ergodic then it would be  isolated from the weak$^*$ closure ${\mathcal K}(\mu)$ of the set of its ergodic components. Since $\mu \in {\mathcal K}(\mu)$,   $\mu$ is ergodic. Therefore,  by definition of ergodicity, $p\omega(x) = \{\mu\}$ for $\mu$-a.e. $x \in S^1$. From the condition $m \ll \mu$ we deduce that $p\omega(x) = \{\mu\}$ for $m$-a.e. $x \in S^1$. This implies, joint with Proposition \ref{PropositionObservable} and Definition    \ref{definitionobservable}, that $\mu$ is the unique SRB-like measure. Now, to end the proof of (b)  it is  only left to check that $\mu \ll m$. Take any Borel set $B \subset S^1$ such that $\mu(B) >0$ and construct the   set $C= \bigcup_{j= 0}^n (f^{-n})(B)$. It satisfies $f^{-1}(C) \subset C$. Since $\mu$ is ergodic and $\mu (C) \geq \mu (B) >0$, we have $\mu(C)= 1$. As $m \ll \mu$ we deduce $m(C) >0$. Therefore $m (f^{-n}(B)) >0$ for some $n \geq 0$. Note that $f^*m \ll   m$, i.e. $m(f^{-1}(B)) = 0$ if $m ( B)= 0$ (this assertion holds because $f \in C^1$ and $f'$ is bounded away from zero). We conclude the $m (B) >0$.
 This shows that $m(B)> 0$ if $\mu(B) > 0$, or in other words $\mu \ll m$, ending the proof. $\ \ \Box$

 \end{prueba}

\vspace{.3cm}

\begin{prueba}
\label{proofCorollary3} {\bf Proof of Corollary \ref{corolario3}} \em

{\em Proof: } From Definition \ref{definitionErgodicSet}    we have $\lim_n (1/n) \sum_{j= 0}^{n-1} \delta_{f^j(x)} = \mu_x \in {\mathcal K}(\mu)$  for $\mu$-almost all the points $x \in S^1$. Since $\mu \ll m$ we have an $m$-positive set of initial states $x \in S^1 $ such that $\{\mu_x\} = p \omega (x) \subset {\mathcal K}(\mu)$. By Definition \ref{definitionobservable} of the set ${\mathcal O}_f$ of SRB-like measures, and after the minimality of ${\mathcal K}(\mu)$ in Definition \ref{definitionErgodicSet}, we deduce that ${\mathcal K}(\mu)  \subset {\mathcal O}_f$.   In other words, if $\nu  \in {\mathcal K}(\mu)$, then $\nu$ is SRB-like.   Applying theorem \ref{maintheorem} and recalling Assertion (\ref{equationRLYequality}), we obtain    ${\mathcal K}(\mu) \subset ES_f$. From  Birkhoff Ergodic Theorem, for any continuous function $\varphi $
$$\int \varphi \ d\mu = \int \lim_{n \rightarrow + \infty} \frac{1}{n} \sum _{j= 0}^{n-1} \varphi \circ f^j \, d \mu = \int \Big (\int \varphi \, d \mu_x \Big ) \, d \mu, $$
where $\mu_x \in {\mathcal K}(\mu)$ is defined by Equality (\ref{equation99}).
    The above integral decomposition  implies that $\mu$ is in the weak$^*$-compact convex hull of ${\mathcal K}(\mu) \subset ES_f$. Since $ES_f$ is weak$^*$-compact and convex (because $f$ is expansive), then $\mu \in ES_f$, as wanted.

    \hfill $\ \ \Box$
\end{prueba}

\end{document}